
\def\LaTeX{%
  \let\Begin\begin
  \let\End\end
  \let\salta\relax
  \let\finqui\relax
  \let\futuro\relax}

\def\UK{\def\our{our}\let\sz s}
\def\USA{\def\our{or}\let\sz z}



\LaTeX

\USA


\salta

\documentclass[twoside,12pt]{article}
\setlength{\textheight}{24cm}
\setlength{\textwidth}{16cm}
\setlength{\oddsidemargin}{2mm}
\setlength{\evensidemargin}{2mm}
\setlength{\topmargin}{-15mm}
\parskip2mm


\usepackage{amsmath}
\usepackage{amsthm}
\usepackage{amssymb}
\usepackage[mathcal]{euscript}
%
%
\usepackage{cite}
\usepackage[usenames,dvipsnames]{color}
%
%
%
\definecolor{darkviolet}{rgb}{0.58, 0.0, 0.83}
\def\gianni{\color{blue}}            
\def\betti{\color{red}}
\def\juerg{\color{Green}}               
\def\pier{\color{darkviolet}}               
\def\new{\color{red}}

%
%

\let\gianni\relax
\let\betti\relax
\let\juerg\relax
\let\pier\relax
\let\new\relax


\bibliographystyle{plain}


%

\finqui

\def\Beq{\Begin{equation}}
\def\Eeq{\End{equation}}
\def\Bsist{\Begin{eqnarray}}
\def\Esist{\End{eqnarray}}

\def\Bthm{\Begin{theorem}}
\def\Ethm{\End{theorem}}

\def\Bex{\Begin{example}\rm}
\def\Eex{\End{example}}

\def\Bcenter{\Begin{center}}
\def\Ecenter{\End{center}}
\let\non\nonumber




\def\step #1 \par{\medskip\noindent{\bf #1.}\quad}

\def\Step #1 \par{\bigskip\leftline{\bf #1}\nobreak\medskip}


\def\Lip{Lip\-schitz}
\def\holder{H\"older}
\def\aand{\quad\hbox{and}\quad}

\def\lhs{left-hand side}
\def\rhs{right-hand side}

\def\wk{well-known}


\def\generaliz{generali\sz}

\def\organiz{organi\sz}

\def\bhv{behavi\our}


\def\multibold #1{\def\arg{#1}%
  \ifx\arg\pto \let\next\relax
  \else
  \def\next{\expandafter
    \def\csname #1#1#1\endcsname{{\bf #1}}%
    \multibold}%
  \fi \next}

\def\pto{.}

\def\multical #1{\def\arg{#1}%
  \ifx\arg\pto \let\next\relax
  \else
  \def\next{\expandafter
    \def\csname cal#1\endcsname{{\cal #1}}%
    \multical}%
  \fi \next}


\def\multimathop #1 {\def\arg{#1}%
  \ifx\arg\pto \let\next\relax
  \else
  \def\next{\expandafter
    \def\csname #1\endcsname{\mathop{\rm #1}\nolimits}%
    \multimathop}%
  \fi \next}

\multibold
qwertyuiopasdfghjklzxcvbnmQWERTYUIOPASDFGHJKLZXCVBNM.

\multical
QWERTYUIOPASDFGHJKLZXCVBNM.

\multimathop
dist div dom meas sign supp .


\def\accorpa #1#2{\eqref{#1}--\eqref{#2}}
\def\Accorpa #1#2 #3 {\gdef #1{\eqref{#2}--\eqref{#3}}%
  \wlog{}\wlog{\string #1 -> #2 - #3}\wlog{}}


\def\separa{\noalign{\allowbreak}}

\def\todx{\mathrel{\scriptstyle\searrow}}

\def\QED{\hfill $\square$}

\def\graffe #1{\mathopen\{#1\mathclose\}}

\def\<#1>{\mathopen\langle #1\mathclose\rangle}
\def\norma #1{\mathopen \| #1\mathclose \|}

\def\normaV #1{\norma{#1}_V}
\def\normaH #1{\norma{#1}_H}

\def\normaVp #1{\norma{#1}_*}

\def\iot {\int_0^t}

\def\iO{\int_\Omega}
\def\intQt{\int_{Q_t}}
\def\intQ{\int_Q}

\def\dt{\partial_t}
\def\dn{\partial_\nu}

\def\cpto{\,\cdot\,}

\def\checkmmode #1{\relax\ifmmode\hbox{#1}\else{#1}\fi}

\def\aeQ{\checkmmode{a.e.\ in~$Q$}}
\def\aet{\checkmmode{a.e.\ in~$(0,T)$}}

\def\aat{\checkmmode{for a.a.~$t\in(0,T)$}}


\def\erre{{\mathbb{R}}}




\def\genspazio #1#2#3#4#5{#1^{#2}(#5,#4;#3)}
\def\spazio #1#2#3{\genspazio {#1}{#2}{#3}T0}

\def\L {\spazio L}
\def\H {\spazio H}

\def\C #1#2{C^{#1}([0,T];#2)}

\def\Vp{V^*}


\def\Lx #1{L^{#1}(\Omega)}
\def\Hx #1{H^{#1}(\Omega)}

\def\Luno{\Lx 1}
\def\Ldue{\Lx 2}

\def\Huno{\Hx 1}
\def\Hdue{\Hx 2}


\def\LQ #1{L^{#1}(Q)}


\let\a\alpha
\let\b\beta

\let\phi\varphi

\let\TeXchi\chi                         
\newbox\chibox
\setbox0 \hbox{\mathsurround0pt $\TeXchi$}
\setbox\chibox \hbox{\raise\dp0 \box 0 }
\def\chi{\copy\chibox}

\let\A\calA
\def\Am{\A^{-1}}
\def\Pi{\widehat\pi}


\def\hatC{\widehat C}


\def\muz{\mu_0}
\def\phiz{\phi_0}
\def\sigmaz{\sigma_0}

\def\Bz{B^0}

\def\muab{\mu_{\a,\b}}
\def\phiab{\phi_{\a,\b}}
\def\sigmaab{\sigma_{\a,\b}}
\def\xiab{\xi_{\a,\b}}
\def\Rab{R_{\a,\b}}
\def\etaab{\eta_{\a,\b}}

\def\phiabp{\phi_{\a,\b'}}

\def\xiabp{\xi_{\a,\b'}}
\def\etaabp{\eta_{\a,\b'}}

\def\mua{\mu_\a}
\def\phia{\phi_\a}
\def\sigmaa{\sigma_\a}
\def\xia{\xi_\a}

\def\etaa{\eta_\a}

\def\mub{\mu_\b}
\def\phib{\phi_\b}
\def\sigmab{\sigma_\b}
\def\xib{\xi_\b}

\def\mul{\overline\mu}
\def\phil{\overline\phi}
\def\sigmal{\overline\sigma}
\def\xil{\overline\xi}
\def\Rl{\overline R}

\def\ustar{u^*}
\def\vstar{v^*}


\def\uz{u_0}


\def\Beta{\widehat{\vphantom t\smash B\mskip2mu}\mskip-1mu}


\Begin{document}


\title{{Asymptotic analyses and error estimates\\
for a Cahn--Hilliard type phase field system \\
\betti{modelling tumor growth}}%
\footnote{{\bf Acknowledgment.}\quad\rm  
{\betti The financial support of the FP7-IDEAS-ERC-StG \#256872
(EntroPhase) is gratefully acknowledged by the authors. The present paper 
also benefits from the support of the MIUR-PRIN Grant 2010A2TFX2 ``Calculus of Variations'' for PC and GG, and the GNAMPA (Gruppo Nazionale per l'Analisi Matematica, la Probabilit\`a e le loro Applicazioni) of INdAM (Istituto Nazionale di Alta Matematica) for PC, GG and ER.}}}

\author{}
\date{}
\maketitle
\Bcenter
\vskip-2cm
{\large\sc Pierluigi Colli$^{(1)}$}\\
{\normalsize e-mail: {\tt pierluigi.colli@unipv.it}}\\[.25cm]
{\large\sc Gianni Gilardi$^{(1)}$}\\
{\normalsize e-mail: {\tt gianni.gilardi@unipv.it}}\\[.25cm]
{\large\sc Elisabetta Rocca$^{(2),\,(3)}$}\\
{\normalsize e-mail: {\tt {\betti rocca@wias-berlin.de and elisabetta.rocca@unimi.it}}}\\[.25cm]
{\large\sc J\"urgen Sprekels$^{(2),\,(4)}$}\\
{\normalsize e-mail: {\tt sprekels@wias-berlin.de}}\\[.45cm]
$^{(1)}$
{\small Dipartimento di Matematica ``F. Casorati'', Universit\`a di Pavia}\\
{\small via Ferrata 1, 27100 Pavia, Italy}\\[.2cm]
$^{(2)}$
{\small Weierstrass Institute for Applied Analysis and Stochastics}\\
{\small Mohrenstrasse 39, 10117 Berlin, Germany}\\[2mm]
$^{(3)}$
{\small Dipartimento di Matematica ``F. Enriques'', Universit\`a di Milano}\\
{\small Via Saldini 50, 20133 Milano, Italy}\\[2mm]
$^{(4)}$
{\small Department of Mathematics,  Humboldt-Universit\"at zu Berlin}\\
{\small {\betti Unter den Linden 6, 10099} Berlin, Germany}\\
[1cm]
\Ecenter
\Begin{abstract}
{\pier This paper is concerned with a phase field system of Cahn--Hilliard type that is related 
to a tumor growth model and consists of three equations in {\gianni terms} of 
the variables order parameter, chemical potential and nutrient concentration. 
This system has been investigated in the recent papers \cite{CGH}
and \cite{CGRS} {\gianni from} the viewpoint of well-posedness, long time \bhv\ 
and asymptotic convergence as two positive viscosity coefficients tend to zero at the same time. 
Here, we  continue the analysis performed in \cite{CGRS} by showing 
two independent sets of results as just one of the coefficents tends to zero, the other remaining fixed. 
We prove convergence results, uniqueness of solutions to the two resulting limit problems, and suitable error estimates.}

\vskip3mm

\End{abstract}

\salta

\pagestyle{myheadings}
\newcommand\testopari{\sc Colli \ --- \ Gilardi \ --- \ Rocca \ --- \ Sprekels}
\newcommand\testodispari{\sc Asymptotic analyses for a Cahn--Hilliard type phase field system}
\markboth{\testodispari}{\testopari}

\finqui


\section{Introduction}
\label{Intro}
\setcounter{equation}{0}

{\pier In this paper, we deal with a system of partial differential equations related to a model for tumor growth {\juerg that was recently proposed in \cite{HZO} (cf.~also \cite{HKNV} and \cite{WZZ}) and further studied analytically} in \cite{CGH,CGRS,FGR}. The system reads} 
\Bsist
  && \a \dt\mu + \dt\phi - \Delta\mu
  = p(\phi) (\sigma - \gamma\mu)
  \label{Iprima}
  \\
  && \mu = \b\dt\phi - \Delta\phi + F'(\phi)
  \label{Iseconda}
  \\
  && \dt\sigma - \Delta\sigma
  = - p(\phi) (\sigma - \gamma\mu)
  \label{Iterza}
\Esist
{\pier and it is complemented} with the boundary and initial conditions
\Bsist
  && \dn\mu = \dn\phi = \dn\sigma = 0 
  \label{Ibc}
  \\
  && \mu(0) = \muz, \quad
  \phi(0) = \phiz
  \aand
  \sigma(0) = \sigmaz \,.
  \label{Icauchy}
\Esist
\Accorpa\Ipbl Iprima Icauchy
Each of the partial differential equations \accorpa{Iprima}{Iterza}
is meant to hold in a three-dimensional bounded domain $\Omega$
endowed with a smooth boundary~$\Gamma$ and for every positive time,
and $\dn$ in \eqref{Ibc} stands for the {\juerg outward} normal derivative on~$\Gamma$.
Moreover, $\a$~and $\b$ are nonnegative parameters, strictly positive in principle,
{\pier and $\gamma$ is a positive} constant.
Furthermore, $p$~is a nonnegative function and $F$ is a nonnegative double-well potential.
Finally, $\muz$, $\uz$ and $\sigmaz$ are given initial data defined in~$\Omega$.

As {\juerg sketched} in~\cite{CGH}
(see also its reference list),
the physical context is that of a tumor-growth model.
The unknown function $\phi$ is an order parameter which is close to two values in the regions of nearly pure phases, 
say, $\phi\simeq 1$ in the tumorous phase and $\phi\simeq -1$ in the healthy cell phase; 
the second unknown $\mu$ is the related chemical potential, specified by \eqref{Iseconda} 
as in the case of the viscous Cahn--Hilliard or Cahn--Hilliard equation, depending on whether 
$\b>0 $ or $\b=0$ (see \cite{CH, EllSt, EllZh});  
the third unknown $\sigma$ stands for the nutrient concentration, 
typically, $\sigma\simeq 1$ in a nutrient-rich extracellular water phase 
and $\sigma\simeq 0$ in a nutrient-poor extracellular water phase. 

In the paper \cite{CGRS},
the asymptotic analysis
as both parameters $\alpha$ and $\beta$ tend to zero at the same time has been performed,
while the {\pier cases} {\juerg when each parameter tends to zero separately {\pier have} been left unanswered}.
{\betti The present paper {\juerg addresses} {\pier these questions; in both cases, the convergence results, 
the uniquenes for the limit problems, as well as the error estimates for the difference of solutions in suitable norms are proved. 
In particular, let us detail here the main mathematical difficulties one encounters} in the two passages to the limit. 

\paragraph{Passage to the limit as $\beta {\pier {} \searrow{} } 0$.}
{\pier The passage to the {\gianni limit} as $\beta$ tends to zero works} 
under quite general assumptions on the proliferation function $p$ and on the 
interaction {\gianni potential~$F$}: 
$p$~is required to be a {\juerg nonnegative}, 
bounded and Lipschitz continuous function, while $F$ is {\juerg the sum of a} convex 
(possibly multivalued)  potential $\Beta$ and a {\pier possibly} nonconvex but smooth part~$\Pi$. 
These assumptions are satisfied in quite a large number of physically meaningful cases like  
the classical double-well potential and the logarithmic potential{\pier ,} defined~by
\Bsist
  \hskip-.8cm && F_{cl}(r) := {\textstyle \frac 14} (r^2-1)^2
  = {\textstyle \frac 14} ((r^2-1)^+)^2 + {\textstyle \frac 14} ((1-r^2)^+)^2
  \quad \hbox{for $r\in\erre$}
  \label{clW}
  \\
  \hskip-.8cm && F_{log}(r) := (1-r)\ln(1-r) + (1+r)\ln(1+r) + \kappa (1 - r^2 )^+
  \quad \hbox{for $|r|<1$}\,,
  \qquad
  \label{logW}
\Esist
where the decomposition $F=\Beta+\Pi$  is written {\juerg explicitly}.
In \eqref{logW}, $\kappa$~is a positive constant
which does or does not provide a double well depending on its value,
and the definition of the logarithmic part of $F_{log}$
is extended by continuity {\juerg to} $\pm1$ and by $+\infty$ outside~$[-1,1]$.
Moreover, another possible choice is the following:
\Beq
  F(r) := I(r) + ((1 - r^2)^+)^2
  \quad \hbox{for $r\in\erre$}\,,
  \label{irrW}
\Eeq
where $I$ is the indicator function of $[-1,1]$,
taking the value $0$ in $[-1,1] $ and $+\infty$ elsewhere.
Regarding the function $p$, in the original model studied in \cite{HZO} 
it was {\pier set as proportional to  $\sqrt{F(\varphi)}$} {\juerg for} $|\varphi|<1$ and zero elsewhere. 
Here we can allow such a behavior as well as more general ones. 
Moreover, we can {\juerg show the uniqueness} of the {\pier solution to the limit problem with $\a>0$ and $\b=0$},  
and an error estimate of the expected order $1/2$ in~$\beta$, 
under a condition of smallness of the fixed {\gianni coefficient~$\alpha$:
the bigger is the Lipschitz constant $\new L$ of the function~$\pi$,
the smaller has to be~$\alpha$}.
{\new About this concern, we construct an example of severe non-uniqueness for the limit problem
(Example~\ref{Nonuniqueness}, where $\a L=1$)}.
 
\paragraph{{\gianni Passage} to the limit as $\alpha{\pier {}\searrow{}} 0$.}  In {\pier this} second case, 
when $\beta$ is kept fixed and $\alpha$ {\juerg tends} to zero, we obtain similar results, but in a 
less general setting. Indeed, in this case, due to the difficulties in estimating the mean value of the 
chemical potential $\mu $ (cf.~\eqref{stimaxi}), we need to assume $F$ to be defined on the 
whole {\juerg of} $\erre$ and to satisfy proper growth conditions (cf.~\eqref{growth}), which are, 
in particular, {\pier fulfilled} {\juerg by} the classical double-well potential \eqref{clW} as well as 
{\juerg by} more general polynomially or exponentially growing potentials. Uniqueness and the 
error estimate can be obtained only under more restrictive conditions on {\pier $p$, which must 
substantially be constant, and $F$, for which it is required a polynomial growth with power four 
(cf.~\eqref{clW}).} This is mainly due to the fact that we need to differentiate in time  equation~
\eqref{Iseconda} in order to get these last results (cf.~Section~\ref{ATOZERO} for more details).  
{\pier However, these results can be compared with \cite[Thms.~2.6 and 2.7]{CGRS} where similar 
(but {\gianni less restirictive}) conditions have to be assumed on~$F$.} 

\vspace{.5cm}

{\pier We think that the methods used in our asymptotic analyses 
could be useful also in other situations. In particular, let us point out that 
in the case of the choice $p\equiv 0$ (admitted by our assumptions 
\eqref{hpp} and \eqref{pconstant})
our system \eqref{Iprima}--\eqref{Iterza} decouples and \eqref{Iprima}--\eqref{Iseconda}
reduces to a well-known phase field system of Caginalp type which can be seen as a 
(doubly) viscous approximation of the Cahn--Hilliard system. 
To this concern, let us quote the paper \cite{Ros1},  where similar asymptotic analyses are carried 
out in the case $p\equiv 0$ with the choice \eqref{clW} for $F$, as both the pameters tend to 0 
or just $\alpha$ goes to $0$ with $\beta >0$ fixed. To the best of our knowledge, we do not know of other investigations  as $
\beta \searrow 0$: this contribution by us seems to be new also in the case $p\equiv 0$. Other examples of rigorous asymptotic analyses with respect to parameters intervening on phase field models can be found, e.g.,  in \cite{CC1, CC2, CGG, CGPS, CS, DKS, Gir, Roc, Ros2, Sch}.} 

\paragraph{Plan of the paper. } Our paper is \organiz ed as follows.}
In the next section, we state our assumptions and results on the mathematical problem.
The last two sections are devoted to the corresponding proofs.


\section{Statement of the problem and results}
\label{STATEMENT}
\setcounter{equation}{0}

In this section, we make precise assumptions and state our results.
As in the Introduction, $\Omega\subset\erre^3$~is the domain where the evolution takes place, and $\Gamma$ is {\juerg its} boundary.
We assume $\Omega$ to be open, bounded and connected, and $\Gamma$ to be smooth.
Moreover, the symbol $\dn$ denotes the outward normal derivative on~$\Gamma$.
Given a final time~$T$, {\pier let} 
\Beq
  Q := \Omega\times(0,T)
  \aand
  \Sigma := \Gamma \times (0,T) .
  \label{defQ}
\Eeq
Moreover, we set for brevity
\Beq
  V := \Huno,
  \quad H := \Ldue
  \aand
  W := \graffe{v\in\Hdue:\ \dn v = 0 \ \hbox{on $\Gamma$}}
  \label{defspazi}
\Eeq
and endow these spaces with their standard norms.
For the norm in a generic Banach space~$X$ (or~a power of~it), we use the symbol~$\norma\cpto_X$
with the following exceptions: we simply write $\norma\cpto_p$ and $\normaVp\cpto$
if $X=\Lx p$ or $X=\LQ p$ for $p\in[1,+\infty]$ and $X=\Vp$, the dual space of~$V$, respectively.
Finally, it is understood that $H\subset\Vp$ as usual, i.e.,
in order that $\<u,v>=\iO uv$ for every $u\in H$ and $v\in V$,
where $\<\cpto,\cpto>$ stands for the duality pairing between $\Vp$ and~$V$.

\bigskip

As far as the structure of the system is concerned, we are given
two constants $\a$ and $\b$ and three functions $p$, $\Beta$ and $\Pi$
satisfying the conditions listed below
\Bsist
  && \a ,\, \b \in (0,1) 
  \label{hpconst}
  \\
  && \hbox{$p:\erre\to\erre$ is nonnegative, bounded and \Lip\ continuous}
  \label{hpp}
  \\
  && \hbox{$\Beta:\erre\to[0,+\infty]$ is convex, proper, lower
      semicontinuous}
    \qquad
  \label{hpBeta}
  \\
  && \hbox{$\Pi\in C^1(\erre)$ is nonnegative and
      $\pi:=\Pi\,'$ is \Lip\ continuous}.
  \label{hpPi}
\Esist
We also define the potential $F:\erre\to[0,+\infty]$
and the graph $B$ in $\erre\times\erre$ by
\Beq
  F := \Beta + \Pi 
  \aand
  B := \partial\Beta 
  \label{defFB}
\Eeq
\Accorpa\HPstruttura hpp defFB
and denote by $D(B)$ and $D(\Beta)$ 
the effective domains of $B$ and~$\Beta$, respectively.
It is well known that $B$
and the operators (denoted with the same symbol~$B$)
induced by $B$ on $L^2$ spaces are maximal monotone
(see, e.g., \cite[Ex.~2.3.4, p.~25]{Brezis}).

We notice that, among many others,
the most important and typical examples of potentials fit our assumptions.
Namely, we can take as $F$ the classical double-well potential
{\betti or} the logarithmic potential defined~{\betti in \eqref{clW} and \eqref{logW}, respectively. 
Another possible choice is the nonsmooth potential \eqref{irrW}.}
In the case of a so irregular potential,
its subdifferential is multi-valued and the precise statement of problem \Ipbl\
has to introduce a selection $\xi$ of~$B(u)$.

As far as the initial data of our problem are concerned,
we assume that
\Beq
  \muz \, , \, \sigmaz \in H , \quad
  \phiz \in V
  \aand
  F(\phiz) \in \Luno\,,
  \label{hpdati}
\Eeq
while the regularity properties we pretend for the solution
are the following:
\Bsist
  && \mu , \sigma \in \H1\Vp \cap \L2V 
  \label{regmusigma}
  \\
  && \phi \in \H1H \cap \L2W 
  \label{regphi}
  \\
  && \xi \in \L2H .
  \label{regxi}
\Esist
\Accorpa\Regsoluz regmusigma regxi
We notice that \accorpa{regmusigma}{regphi} imply
$\mu,\sigma\in\C0H$ and $\phi\in\C0V$.
At this point, we consider the problem 
of finding a quadruplet $(\mu,\phi,\sigma,\xi)$
with the above regularity 
in order that $(\mu,\phi,\sigma,\xi)$ and the related function
\Beq
  R = p(\phi) (\sigma - \mu)
  \label{defR}
\Eeq
satisfy the system
\Bsist
  \hskip-1cm
  && \a \< \dt\mu , v > + \iO \dt\phi \, v + \iO \nabla\mu \cdot \nabla v
  = \iO R v 
  \non
  \\
  && \quad \hbox{for every $v\in V$, \ \aet} 
  \label{prima}
  \\ [3pt]
  && \mu = \b\dt\phi - \Delta\phi + \xi + \pi(\phi)
  \aand
  \xi \in B(\phi)
  \quad \aeQ
  \label{seconda}
  \\
  && \< \dt\sigma , v > + \iO \nabla\sigma \cdot \nabla v
  = - \iO R v 
  \non
  \\
  && \quad \hbox{for every $v\in V$, \ \aet} 
  \label{terza}
  \\ [3pt]
  && \mu(0) = \muz, \quad
  \phi(0) = \phiz
  \aand
  \sigma(0) = \sigmaz \,.
  \label{cauchy}
\Esist
\Accorpa\Pbl defR cauchy
\Accorpa\Tuttopbl regmusigma cauchy
This is a weak formulation of the boundary value problem \Ipbl\ described in the Introduction.
The homogeneous Neumann boundary condition for $\phi$ 
is contained in~\eqref{regphi} (see \eqref{defspazi} for the definition of~$W$),
while the analogous ones for $\mu$ and $\sigma$
are meant in a \generaliz ed sense through the variational equations \eqref{prima} and~\eqref{terza}.
We notice once and for all that 
{\juerg the addition of \eqref{prima} and~\eqref{terza}} yields
\Beq
  \< \dt \bigl( \a \mu + \phi + \sigma \bigr) , v >
  + \iO \nabla (\mu + \sigma) \cdot \nabla v
  = 0
  \label{unopiutre}
\Eeq
for every $v\in V$, \aet.

Well-posedness for the above problem is ensured by \cite[Thm.~2.2]{CGRS},
which {\juerg states} 

\Bthm
\label{Wellposedness}
Assume \HPstruttura\ and \eqref{hpdati}.
Then, for every $\a,\b\in(0,1)$, there exists a unique quadruplet $(\mu,\phi,\sigma,\xi)$
satisfying \Regsoluz\ and solving problem \Pbl.
\Ethm

{\juerg In the} same paper~\cite{CGRS}, the authors study the asymptotic analysis of the above problem
as both the parameters $\a$ and $\b$ tend to zero {\juerg simultaneously}
and prove an error estimate for the difference of the solution to problem \Pbl\
and the one of the expected limit problem under further assumptions on the potential~$F$
(see \cite[Thms.~2.5 and~2.6]{CGRS}).
Namely, it is assumed that $F$ is everywhere defined and smooth
and satisfies some growth condition.
In particular, the classical double-well potential \eqref{clW} is allowed.

In the present paper, we discuss the analogous problems obtained
by {\juerg letting} just one of the parameters tend to zero
{\juerg while} keeping the other fixed.
The results that deal with the possible cases are presented at once.

In the first {\juerg case, we keep} $\a$ fixed and let $\b$ tend to zero.
The limit problem one expects is the following:
\Bsist
  & \< \dt(\a\mu+\phi) , v > + \iO \nabla\mu \cdot \nabla v
  = \iO R \, v
  & \quad \hbox{$\forall\,v\in V$, \ \aet}
  \label{primaa}
  \\
  & \mu = - \Delta\phi + \xi + \pi(\phi)
  \aand
  \xi \in B(\phi)
  & \quad \aeQ
  \label{secondaa}
  \\
  & \< \dt\sigma , v > + \iO \nabla\sigma \cdot \nabla v
  = - \iO R \, v
  & \quad \hbox{$\forall\,v\in V$, \ \aet}
  \label{terzaa}
  \\
  & (\a\mu+\phi)(0) = \a\muz + \phiz
  \aand
  \sigma(0) = \sigmaz
  & \quad \hbox{in $\Omega$}\,,
  \label{cauchya}
\Esist
\Accorpa\Pbla primaa cauchya
where $R$ is defined by~\eqref{defR}.
For its solution, we require that
\Bsist
  && \mu \in \L\infty H \cap \L2V
  \label{regmua}
  \\
  && \phi \in \L\infty V \cap \L2W
  \label{regphia}
  \\
  && \a\mu + \phi \in \H1\Vp
  \label{regamupiuphi}
  \\
  && \sigma \in \H1\Vp \cap \L2V 
  \label{regsigmaa}
  \\
  && \xi \in \L2H \,.
  \label{regxia}
\Esist
\Accorpa\Regsoluza regmua regxia
Our result on the asymptotic \bhv\ as $\b\todx0$ holds
{\pier under the assumption that $\alpha $ is sufficiently small}.

\Bthm
\label{AsymptoticsB}
Assume \HPstruttura\ on the structure and \eqref{hpdati} on the initial data. {\pier Then, there exists $\a_0\in(0,1)$ such that, for $\a\in(0,\a_0)$ and  $\b\in(0,1)$, the unique solution $(\muab,\phiab,\sigmaab,\xiab)$ to problem \Pbl, with the regularity \Regsoluz, satisfies} 
\Bsist
  & \muab \to \mua 
  & \hbox{weakly {\pier star} in $\L\infty H\cap\L2V$}
  \label{bconvmu}
  \\
  & \phiab \to \phia
  & \hbox{weakly star in $\L\infty V\cap\L2W$} \qquad\quad
  \label{bconvphi}
  \\
  & \sigmaab \to \sigmaa
  & \hbox{weakly in $\H1\Vp\cap\L2V$}
  \label{bconvsigma}
  \\
  & \dt (\a\muab+\phiab) \to \dt (\a\mua+\phia)
  & \hbox{weakly in $\L2\Vp$}
  \label{bconvdt}
  \\
  & {\pier \beta \phiab \to 0}
  & {\pier \hbox{strongly in $\H1H\cap\L2W$}}
  \label{pier2}
  \\
  & \xiab \to \xia
  & \hbox{weakly in $\L2H$}
  \label{bconvxi}
\Esist
as $\b$ tends to zero, at least for a subsequence,
and every limiting quadruplet $(\mua,\phia,\sigmaa,\xia)$
solves problem \Pbla.
In particular, problem \Pbla\ has at least a solution satisfying \Regsoluza.
\Ethm

{\pier Moreover,} we can prove uniqueness for the limit problem and an error estimate {\pier under an additional restriction on~$\a$.}

\Bthm
\label{Btozero}
Assume \HPstruttura\ and \eqref{hpdati}.
Then, there exists {\pier $\a_{00}\in(0,\a_0)$ such that, for $\a\in(0,\a_{00})$,} the following conclusions hold:

\noindent
$i)$\enskip the solution $(\mua,\phia,\sigmaa,\xia)$ to problem \Pbla\ satisfying \Regsoluza\ is unique;

\noindent
$ii)$\enskip if $(\muab,\phiab,\sigmaab,\xiab)$ is the unique solution
to problem \Pbl\ with the regularity specified by \Regsoluz\ for $\b\in(0,1)$,
the estimate
\Bsist
  && \norma{\muab-\mua}_{\L2H}
  + \norma{\phiab-\phia}_{\L2V}
  + \norma{\sigmaab-\sigmaa}_{\L\infty H\cap\L2V}
  \qquad
  \non
  \\
  && \quad {}
  + \norma{(\a\muab+\phiab+\sigmaab)-(\a\mua+\phia+\sigmaa)}_{\L\infty\Vp}
  \leq C_\a \, \b^{1/2} 
  \label{errorbtozero}
\Esist
holds true for $\b\in(0,1)$ with a constant $C_\a$ that depends only on~$\a$, $\Omega$, $T$, 
the structure of the system,
and the norms of the initial data related to assumptions~\eqref{hpdati}.
\Ethm

{\new
\Bex
\label{Nonuniqueness}
As said in the Introduction, we can construct an example of severe non-uniqueness
for the limit problem.
We take $p=0$, so that the third equation is decoupled. 
Moreover, we take the indicator function of $[-1,1]$ as~$\Beta$
and, given $L>0$, we choose $\Pi$ smooth, nonnegative and such~that
\Beq
  |\pi'(r)| \leq L
  \quad \hbox{for every $r\in\erre $}
  \aand
  \pi(r) = - Lr
  \quad \hbox{for $r\in [-1,1]$}.
  \non
\Eeq
Finally, we take 
\Beq
  \phiz = \muz = \sigmaz = 0 .
  \non
\Eeq
Then, for any function $\psi\in L^\infty(0,T)$ satisfying $|\psi(t)|\leq1$ \aat, the definition
\Beq
  \mu(x,t) := - L \psi(t) , \quad
  \phi(x,t) := \psi(t) , \quad
  \sigma(x,t) := 0
  \aand
  \xi(x,t) := 0
  \non
\Eeq
provides a solution if $\a L=1$.
Indeed, the initial and boundary conditions are trivially satisfied
as well as the equations, since we have 
\Beq
  \alpha\mu + \phi
  = - \a L \psi + \psi
  = 0 , \quad
  \mu = - L \psi = - L \phi = \pi(\phi)
  \aand
  \xi = 0 \in B(\psi) = B(\phi) .
  \non
\Eeq
\Eex
}

In the second case, we keep $\b$ fixed
and study the asymptotics with respect to the parameter $\a$
and the corresponding expected limit problem, namely
\Bsist
  & \< \dt\phi , v > + \iO \nabla\mu \cdot \nabla v
  = \iO R \, v
  & \quad \hbox{$\forall\,v\in V$, \ \aet}
  \label{primab}
  \\
  & \mu = \b\dt\phi - \Delta\phi + \xi + \pi(\phi)
  \aand
  \xi \in B(\phi)
  & \quad \aeQ
  \label{secondab}
  \\
  & \< \dt\sigma , v > + \iO \nabla\sigma \cdot \nabla v
  = - \iO R \, v
  & \quad \hbox{$\forall\,v\in V$, \ \aet}
  \label{terzab}
  \\
  & \phi(0) = \phiz
  \aand
  \sigma(0) = \sigmaz
  & \quad \hbox{in $\Omega$}\,,
  \label{cauchyb}
\Esist
\Accorpa\Pblb primab cauchyb
where $R$ is defined by~\eqref{defR}.
For its solution, we require the following regularity:
\Bsist
  && \mu \in \L2V
  \label{regmub}
  \\
  && \phi \in \H1H \cap \L2W
  \label{regphib}
  \\
  && \sigma \in \H1\Vp \cap \L2V
  \label{regsigmab}
  \\
  && \xi \in \L2H \,.
  \label{regxib}
\Esist
\Accorpa\Regsoluzb regmub regxib
Our results hold in a less general setting.
We need a first restriction on the potential~$F$ 
in order to prove an asymptotic result
and the existence of a solution to the limit problem.
Namely, we also assume~that
\Beq
  D(\Beta) = \erre
  \aand
  |\Bz(r)| \leq C_B \bigl( \Beta(r) + 1 \bigr)
  \quad \hbox{for every $r\in\erre$}\,,
  \label{growth}
\Eeq
where $\Bz(r)$~is the element of $B(r)$ having minimum modulus
and $C_B$ is a given constant.
Notice that the classical potential \eqref{clW}
and similar potentials with {\juerg polynomial or exponential growth satisfy this} 
 assumption.
Moreover, \eqref{growth}~still allows $B$ to be multi-valued.
We~have {\juerg the following result}:

\Bthm
\label{AsymptoticsA}
Assume \HPstruttura\ and \eqref{growth} on the structure and \eqref{hpdati} on the initial data.
Moreover, for $\a,\b\in(0,1)$, let $(\muab,\phiab,\sigmaab,\xiab)$ 
be the unique to problem \Pbl\ satisfying \Regsoluz.
Then {\juerg we~have that}
\Bsist
  & \muab \to \mub
  & \hbox{weakly in $\L2V$}
  \label{aconvmu}
  \\
  & \phiab \to \phib
  & \hbox{{\pier weakly in} $\H1H\cap\L2W$} \qquad
  \label{aconvphi}
  \\
  & \sigmaab \to \sigmab
  & \hbox{weakly in $\H1\Vp\cap\L2V$}
  \label{aconvsigma}
  \\
  & \xiab \to \xib
  & \hbox{weakly in $\L2H$}  
  \label{aconvxi}
  \\
  & {\pier \a\muab \to 0}
  & {\pier \hbox{weakly in $\H1\Vp$ and strongly in $\L2V$}}\qquad
  \label{aconvdt}
\Esist
as $\a$ tends to zero, at least for a subsequence,
and every limiting quadruplet $(\mub,\phib,\sigmab,\xib)$
solves problem \Pblb.
In particular, the limit problem \Pblb\ has at least a solution satisfying \Regsoluzb.
\Ethm

Both {\juerg the} uniqueness for the limit problem and 
{\juerg the} error estimates look difficult to prove.
We can overcome such a difficulty only in a particular case.
Moreover, we have to {\pier take $p$ constant} and make stronger conditions on the potential~$F$
(which, however, still {\juerg allow} the choice of the classical potential~\eqref{clW})
and~slightly reinforce those on the initial data.
Namely, we also assume that
\Bsist
  && \hbox{$p$ is a {\pier nonnegative} constant}
  \label{pconstant}
  \\[0.1cm]
  && \hbox{$F$ is a $C^2$ function on $\erre$ satisfying } \, |F''(r)| 
  \leq C (r^2+1)
  \non
  \\
  &&
  \qquad \hbox{for every $r\in\erre$ {\pier and for some constant $C>0$.}}
  \label{classico}
\Esist
We observe that {\pier in this setting 
the last condition in \eqref{hpdati} is a consequence of 
$\phiz \in V$ and \eqref{classico}.}
Indeed, these assumptions ensure that $F(r)=O(r^4)$ as $|r|$ tends to~$+\infty$ and that $\phiz\in\Lx4$, thanks to the Sobolev inequality.

\Bthm
\label{Atozero}
Assume \HPstruttura, \eqref{growth} and \accorpa{pconstant}{classico} on the structure,
and {\pier \eqref{hpdati}} on the initial data,
and let $\b\in(0,1)$.
Then, the following conclusions hold {\juerg true}:

\noindent
$i)$\enskip the solution $(\mub,\phib,\sigmab,\xib)$ to problem \Pblb\ satisfying \Regsoluzb\ is unique;

\noindent
$ii)$\enskip if $(\muab,\phiab,\sigmaab,\xiab)$ is the unique solution
to problem \Pbl\ with the regularity specified by \Regsoluz\ for $\a\in(0,1)$,
the estimate
\Bsist
  && \norma{\muab-\mub}_{\L2V}
  + \norma{\phiab-\phib}_{\L\infty H\cap\L2V}
  \non
  \\
  && \quad {}
  + \norma{\sigmaab-\sigmab}_{\L\infty H\cap\L2V}
  \leq C_\b \, \a^{1/2} 
  \label{erroratozero}
\Esist
holds true for $\a\in(0,1)$ with a constant $C_\b$ that depends only on~$\b$, $\Omega$, $T$, 
the structure of the system,
and the norms of the initial data related to assumptions~{\pier\eqref{hpdati}}.
\Ethm

Now, we list some facts.
We repeatedly make use of the notation
\Beq
  Q_t := \Omega \times (0,t)
  \quad \hbox{for $t\in[0,T]$}
  \label{defQt}
\Eeq
and of \wk\ inequalities, namely,  
the elementary Young inequality
\Beq
  ab \leq \delta a^2 + \frac 1{4\delta} \, b^2
  \quad \hbox{for every $a,b\geq 0$ and $\delta>0$}\,,
  \label{young}
\Eeq
the \holder\ inequality, and its consequences.
Moreover, as $\Omega$ is bounded and smooth,
we can owe to the Poincar\'e and Sobolev--type inequalities,
namely, 
\Bsist
  && \normaV v \leq C \Bigl( \normaH{\nabla v} + \bigl| \textstyle\iO v \bigr| \Bigr)
  \quad \hbox{for every $v\in V$} 
  \label{poincare}
  \\
  && V \subset \Lx q
  \aand
  \norma v_q \leq C \normaV v
  \quad \hbox{for every $v\in V$ and $1\leq q\leq 6$}
  \label{sobolev}
  \\
  && \Lx q \subset \Vp
  \aand
  \normaVp v \leq C \norma v_q
  \quad \hbox{for every $v\in\Lx q$ and $q\geq 6/5$} \,.
  \qquad
  \label{dualsobolev}
\Esist
In \accorpa{poincare}{dualsobolev}, $C$ only depends on~$\Omega$.
Furthermore, we denote by $\A$ the Riesz isomorphism from $V$ onto $\Vp$
associated to the standard inner product of~$V$, i.e.,
$\A:V\to\Vp$ is defined~by
\Beq
  \< \A u, v > := (u,v)_V
  = \iO \bigl( \nabla u \cdot \nabla v + uv \bigr)
  \quad \hbox{for $u,v\in V$}.
  \label{defA}
\Eeq
We notice that $\A u=-\Delta u+u$ if $u\in W$.
We also remark that
\Bsist
  && \< \A u , \Am \vstar >
  = \< \vstar , u >
  \quad \hbox{for every $u\in V$ and $\vstar\in\Vp$}
  \label{propAuno}
  \\
  && \< \vstar , \Am \ustar >
  = \< \ustar , \Am \vstar >
  = (\ustar,\vstar)_*
  \quad \hbox{for every $\ustar,\vstar\in\Vp$}\,,
  \label{propAdue}
\Esist
where $(\cpto,\cpto)_*$ {\juerg denotes} the dual scalar product in $\Vp$
associated to the standard one in~$V$,
and recall that $\<\vstar,u>=\iO\vstar u$ if $\vstar\in H$.
As a consequence of~\eqref{propAdue}, we~have
\Beq
  \frac d{dt} \, \normaVp\vstar^2
  = 2 \< \dt\vstar , \Am \vstar >
  \quad \hbox{for every $\vstar\in\H1\Vp$} .
  \label{propAtre}
\Eeq
Finally, in order to simplify the notation, 
we follow a general rule in performing our a~priori estimates.
The small-case italic $c$ without any subscript stands for different constants,
which may only depend on~$\Omega$, $T$, 
the shape of the nonlinearities 
and the norms of the initial data related to assumptions at hand.
A~notation like~$c_\delta$ signals a constant that {\juerg also depends}
on the parameter~$\delta$.
We point out that $c$ and $c_\delta$ do not depend on $\a$ and~$\b$
and that their meaning might change from line to line and even 
{\juerg within} the same chain of inequalities.
On the contrary, {\juerg constants that are later referred to} are always denoted by different symbols,
e.g., by a capital letter.

The starting point for the proofs {\juerg given} in the next sections
is the following general result (see \cite[Thm.~2.3]{CGRS}):

\Bthm
\label{GenEst}
Assume \HPstruttura\ and \eqref{hpdati}.
Then, for some constant $\hatC$ {\juerg which only depends on} 
$\Omega$, $T$ and the shape of the nonlinearities,
the following {\juerg holds true: for} every $\a,\b\in(0,1)$,
the solution $(\mu,\phi,\sigma,\xi)$
to problem \Pbl\ with the regularity specified by~\Regsoluz\ satisfies
\Bsist
  && \a^{1/2} \norma\mu_{\L\infty H}
  + \norma{\nabla\mu}_{\L2H}
  \non
  \\
  && \quad {}
  + \b^{1/2} \norma{\dt\phi}_{\L2H}
  + \norma\phi_{\L\infty V}
  + \norma{F(\phi)}_{\L\infty\Luno}^{1/2}
  \non
  \\
  && \quad {}
  + \norma{\dt(\a\mu+\phi)}_{\L2\Vp}
  + \norma\sigma_{\H1\Vp\cap\L\infty H\cap\L2V}
  \non
  \\
  && \leq \hatC \, \bigl(
    \a^{1/2} \normaH\muz + \normaV\phiz + \norma{F(\phiz)}_{\Luno}^{1/2} + \normaH\sigmaz 
  \bigr)
  \label{genest}
\Esist
as well~as
\Bsist
  && \norma\mu_{\L2V}
  + \norma\phi_{\L2W}
  + \norma\xi_{\L2H}
  \non
  \\
  && \leq \hatC \, \bigl(
    \a^{1/2} \normaH\muz + \normaV\phiz + \norma{F(\phiz)}_{\Luno}^{1/2} + \normaH\sigmaz + \norma\mu_{\L2H} + 1
  \bigr) .
  \qquad
  \label{genestbis}
\Esist
\Ethm


\section{{\pier Proofs} of Theorems~\ref{AsymptoticsB} and~\ref{Btozero}}
\label{BTOZERO}
\setcounter{equation}{0}
 
It is understood that the assumptions \HPstruttura\ and \eqref{hpdati} are in force.

\step
Proof of {\betti Theorem}~\ref{AsymptoticsB} 

We start from {\pier the uniform estimates stated in} Theorem~\ref{GenEst}.
As $\a$ is fixed, \eqref{genest}~ensures that 
even $\norma\muab_{\L\infty H}$ is bounded,
so that \eqref{genestbis} provides a bound for its \lhs.
Therefore the convergence{\betti s} specified in {\pier \accorpa{bconvmu}{bconvxi}} hold for a subsequence.
Moreover, {\pier as $\etaab:= \a \muab +\phiab $ converges to $\etaa:= \a \mua +\phia $ weakly in $\H1\Vp \cap \L2V$, by applying the Aubin-Lions lemma (see, e.g., \cite[Thm.~5.1, p.~58]{Lions}) we deduce that 
\Beq
  \etaab \to \etaa
  \quad \hbox{strongly in $\L2H$ as $\beta \searrow  0$.}
  \label{pier3}
\Eeq
Now, we prove that \eqref{pier3} implies that both
\Beq
  \muab \to \mua \quad \hbox{and} \quad  \phiab \to \phia 
  \quad \hbox{strongly in $\L2H$}
  \label{pier4}
\Eeq
provided $\a$ is small enough. 
We show \eqref{pier4} by using a Cauchy sequence argument: 
{\gianni we write \eqref{seconda} for the solutions corresponding to $\beta$ and~$\beta'$
and take the difference;
then, we multiply by~$\a$, sum $\phiab - \phiabp$ to both sides, rearrange and finally test by $\phiab - \phiabp$}. 
We obtain
{\gianni
\begin{align}
  & \iO \bigl( (\etaab -\etaabp) - 
  {\pier \a(\b\,\dt\phiab - \b'\,\dt\phiabp)} \bigr) \, (\phiab -\phiabp)
  \non
  \\
  & {} = \normaH{\phiab - \phiabp}^2 
  + \a \iO |\nabla (\phiab - \phiabp)|^2
  \non
  \\
  & \quad {} + \a \iO (\xiab -\xiabp)(\phiab - \phiabp) 
  + \a \iO (\pi(\phiab) - \pi(\phiabp))(\phiab - \phiabp)
  \label{pier5}
\end{align}
a.e.\ in $(0,T)$}. 
Then we integrate over $(0,T)$ and observe that the resulting \lhs\ tends to zero as  
$\beta , \beta' \searrow  0$, due to the 
{\gianni strong convergence given in \eqref{pier3} and \eqref{pier2}
coupled with the weak convergence \eqref{bconvphi}}. 
The term $\a \intQ (\xiab -\xiabp)(\phiab - \phiabp)$ is nonegative thanks to the monotonicity of~$B$; 
the last term can be treated using the Lipschitz continuity of~$\pi$, namely, 
$$ \a \intQ (\pi(\phiab) - \pi(\phiabp))(\phiab - \phiabp) \geq - \a L  \norma{\phiab - \phiabp}^2_2 ,$$
where $L$ denotes a Lipschitz constant for~$\pi$. 
Hence, from \eqref{pier5} and the subsequent remarks we deduce that $\{\phiab\}$ is a Cauchy sequence in $L^2(Q)\equiv \L2H $ 
provided  $\alpha L<1$, so that \eqref{bconvphi} entails that $\phiab$ converges to {\gianni $\phia$} strongly in $\L2H$
and \eqref{pier4} is completely proved.}

{\pier From \eqref{pier4} it} follows that $\pi(\phiab)$ and $p(\phiab)$
converge to $\pi(\phia)$ and $p(\phia)$, respectively, strongly in $\L2H$.
Therefore, we can identify the limits of the nonlinear terms $\xiab$ and~$\Rab$.
For the former, we can apply, e.g., \cite[Cor.~2.4, p.~41]{Barbu}.
For the latter, we note that $\Rab$ converges to $p(\phia)(\sigmaa-\mua)$ 
{\pier strongly in~$\LQ1$ since (cf.~\eqref{bconvsigma}) we also have 
a strong convergence of $\sigmaab $ to $\sigmaa$  in $L^2(Q)$.}
At this point, we can write the integrated--in--time version of problem \accorpa{prima}{terza}
for the approximating solution
with time dependent test functions and take the limit as $\a$ tends to zero.
We obtain the analogous systems for $(\mua,\phia,\sigmaa,\xia)$,
and this implies {\betti  \accorpa{primaa}{cauchya}} for such a quadruplet.\QED

\step
{\pier Proof  of} Theorem~\ref{Btozero}

{\pier In order to show Theorem~\ref{Btozero},
we} do not follow the order of the statement.
Indeed, assume for a while that the part~$ii)$ 
has been proved for every solution to problem \Pbla\
(with a constant $C_\a$ that might depend on the solution we are considering),
provided that $\a$ is small.
From this, we derive the uniqueness part~$i)$.
Indeed, let $(\mu_i,\phi_i,\sigma_i,\xi_i)$, $i=1,2$, be two solutions.
If $\a$ is small, inequality \eqref{errorbtozero}
holds for both of them with a common constant~$C_\a$.
Hence, we immediately derive,~e.g.,
\Bsist
  && \norma{\mu_1-\mu_2}_{\L2H}
  + \norma{\phi_2-\phi_2}_{\L2V}
  \non
  \\
  && \quad {}
  + \norma{\sigma_1-\sigma_2}_{\L2V}
  \leq 2 C_\a \, \b^{1/2} 
  \quad \hbox{for every $\b\in(0,1)$}.
  \non
\Esist
This implies that $\mu_1=\mu_2$, $\phi_1=\phi_2$ and $\sigma_1=\sigma_2$.
Then, $\xi_1=\xi_2$ by comparison in \eqref{secondaa}.

{\pier Now,} we {\pier show that the error estimate {\betti \eqref{errorbtozero}}
holds for every solution $(\mua,\phia,\sigmaa,\xia)$ to the} problem \Pbla\
(with a constant $C_\a$ that might depend on the solution we are considering).
To~this end, we present equations 
\accorpa{seconda}{terza}, \accorpa{secondaa}{terzaa},
\eqref{unopiutre} and its analogue 
obtained by {\pier summing \eqref{primaa} and \eqref{terzaa}} 
in a slightly different form.
Namely, in each equation, we add the same function to both sides
and make the Riesz ismorphism~$\A$ appear (see \eqref{defA}).
In doing {\juerg so}, we owe to the regularity conditions \Regsoluz\ and \Regsoluzb. 
If $(\mul,\phil,\sigmal,\xil)$ stands for the solution to the limit problem, we~have
\Bsist
  && \dt(\a\muab + \phiab + \sigmaab) + \A(\muab + \sigmaab)
  = \muab + \sigmaab
  \label{unopiutreA}
  \\
  && \muab = \b\dt\phiab + \A\phiab + \xiab + \pi(\phiab) - \phiab
  \label{secondaA}
  \\
  && \dt\sigmaab + \A\sigmaab = -\Rab + \sigmaab
  \label{terzaA}
  \\
  && \dt(\a\mul + \phil + \sigmal) + \A(\mul + \sigmal)
  = \mul + \sigmal
  \label{unopiutrelA}
  \\
  && \mul = \A\phil + \xil + \pi(\phil) - \phil
  \label{secondalA}
  \\
  && \dt\sigmal + \A\sigmal = - \Rl + \sigmal\,,
  \label{terzalA}
\Esist
where $\Rab$ and $\Rl$ are defined by~\eqref{defR} according
to the equations we are considering.
All these equations are meant in the framework of the Hilbert triplet $(V,H,\Vp)$, \aet.
However, the explicite versions of \eqref{secondaA} and \eqref{secondalA}
also hold \aeQ.
Moreover, we have to add the conditions $\xiab\in B(\phiab)$ and $\xil\in B(\phil)$ \aeQ\
as well as the initial conditions
\eqref{cauchy} and~\eqref{cauchya}, respectively.
Now, we take the differences between 
\accorpa{unopiutreA}{terzaA} and \accorpa{unopiutrelA}{terzalA}
and~have
\Bsist
  && \dt(\a\mu + \phi + \sigma) + \A(\mu + \sigma)
  = \mu + \sigma
  \non
  \\
  && \mu = \b\dt\phiab + \A\phi + \xiab -\xil + \pi(\phiab) - \pi(\phil) - \phi
  \non
  \\
  && \dt\sigma + \A\sigma = - (\Rab - \Rl) + \sigma\,,
  \non
\Esist
where we have set, for convenience,
\Beq
  \mu := \muab - \mul , \quad
  \phi := \phiab - \phil 
  \aand
  \sigma := \sigmaab - \sigmal .
  \non
\Eeq
At this point, we {\juerg write} these equations at time $s\in(0,T)$
and test them~by
\Beq
  \Am (\a\mu + \phi + \sigma)(s) , \quad
  -\phi(s)
  \aand
  \sigma(s)\,,
  \non
\Eeq
respectively. 
Next, we sum up and integrate over $(0,t)$ with respect to~$s$,
for an arbitrary $t\in(0,T)$.
Then, we {\pier recall and use  \accorpa{defA}{propAtre}; by} 
rearranging a little and omitting the evaluation at $s$ inside integrals for brevity,
we~obtain {\juerg that}
\Bsist
  && \frac 12 \, \normaVp{(\a\mu + \phi + \sigma)(t)}^2
  + \intQt (\mu + \sigma) (\a\mu + \phi + \sigma) 
  {\pier {}+ \iot \normaV{\phi}^2 \, ds}
  \non
  \\
  && \quad {\pier {}
  + \intQt ( \xiab - \xil) \phi
  - \intQt \mu \phi 
  + \frac 12 \iO |\sigma(t)|^2
  + \iot \normaV\sigma^2 \, ds}
  \non
  \\
  \separa
  && = \iot \bigl( \mu + \sigma , \a\mu + \phi + \sigma \bigr)_* \, ds 
  - \b \intQt \dt\phiab \, \phi 
  + \intQt \bigl\{ \phi - \bigl( \pi(\phiab) - \pi(\phil) \bigr) \bigr \} \phi 
  \qquad\quad
  \non
  \\
  && \quad {}
  - \intQt (\Rab - \Rl) \sigma 
  + \intQt |\sigma|^2 \,.
  \label{pererroreA}
\Esist
The terms on the \lhs\ {\juerg not having a definite sign
are treated simultaneously in the following way:}
\Bsist
  && \intQt (\mu + \sigma) (\a\mu + \phi + \sigma)
  - \intQt \mu \phi
  = \intQt \bigl( \a |\mu|^2 + \mu\sigma + \sigma(\a\mu + \phi + \sigma) \bigr) 
  \non
  \\
  && \geq \a \intQt |\mu|^2
  + \intQt \mu\sigma
  - \iot \normaV\sigma \, \normaVp{\a\mu + \phi + \sigma} \, ds
  \non 
  \non
  \\
  && \geq \frac \a 2 \intQt |\mu|^2
   - c \intQt |\sigma|^2 
   - \delta \iot \normaV\sigma^2 \, ds
   - c_\delta \iot \normaVp{\a\mu + \phi + \sigma}^2 \, ds .
   \non
\Esist
Now, we deal with the \rhs\ of \eqref{pererroreA}.
For the first integral, we observe that
the norm of the embedding $H\subset\Vp$ is~$1$
(since the norms of $V$ and $H$ are the standard ones),
and {\juerg we have the estimate}
\Bsist
  && \iot \bigl( \mu + \sigma , \a\mu + \phi + \sigma \bigr)_* \, ds
  \leq  {\pier \frac \a 8 }\iot \bigl(\normaVp\mu^2 + \normaVp\sigma^2 \bigr) \, ds
  + {\pier c} \iot \normaVp{\a\mu + \phi + \sigma}^2 \, ds 
  \non
  \\
  && \leq {\pier \frac \a 8 }\intQt |\mu|^2
  +  {\pier \frac \a 8 } \intQt |{\betti \sigma}|^2
  + {\pier c } \iot \normaVp{\a\mu + \phi + \sigma}^2 \, ds \,.  
  \non
\Esist
Next, we~have {\juerg that}
\Beq
  - \b \intQt \dt\phiab \, \phi
  \leq \delta \intQt |\phi|^2
  + c_\delta \, \b^2 \intQt |\dt\phiab|^2
  \leq \delta \iot \normaV\phi^2 \, ds
  + c_\delta \beta \, ,
  \non
\Eeq
the last inquality being due to~\eqref{genest} for $\dt\phiab$.
For the next term, we use the \Lip\ continuity of $\pi$
and {\pier still denote by $L$ the \Lip\ constant of~$\pi$.}
Then, we {\juerg obtain that}
\Bsist
  && \intQt \bigl\{ \phi - \bigl( \pi(\phiab) - \pi(\phil) \bigr) \bigr \} \phi
  \leq {\pier (1+L)}  \intQt |\phi|^2
  \non
  \\
  \separa
  && = {\pier (1+L) \left( \intQt \phi \bigl( \a\mu + \phi + \sigma \bigr)
  - \a \intQt \phi \mu  - \intQt \phi \sigma \right) }
  \non
  \\
  && \leq \delta \iot \normaV\phi^2
  + c_\delta \iot \normaVp{\a\mu + \phi + \sigma}^2 \, ds 
  \non
  \\
  && \quad {}
  + {\pier \frac \a 8} \intQt |\mu|^2
  + {\pier (1+L)}^2 \a \iot \normaV\phi^2 \, ds
  + \delta \iot \normaV\phi^2 \, ds
  + c_\delta \intQt |\sigma|^2
  \non
\Esist
and we assume at once $\a$ to be small in order that ${\pier (1+L)}^2\a<1$.
Finally, {\pier there is one more term to treat} on the \rhs\ of \eqref{pererroreA}, 
{\pier namely,}
\Bsist
  && - \intQt (\Rab - \Rl) \sigma
  = - \intQt \bigl\{ p(\phiab) (\sigmaab - \muab) - p(\phil) (\sigmal - \mul) \bigr\} \sigma
  \non
  \\
  && \leq \intQt |p(\phiab) - p(\phil)| \, |\sigmaab - \muab| \, |\sigma|
  + \intQt |p({\betti\phil})| \, |\sigma - \mu| \, |\sigma|
  \non
  \\
  && \leq c \intQt |\phi| \, |\sigmaab - \muab| \, |\sigma|
  + c \intQt |\sigma - \mu| \, |\sigma|\, .
  \non
\Esist
{\juerg {\pier Note that} the last inequality holds true} since $p$ is \Lip\ continuous and bounded.
On the other hand, we can use the \holder\ and Sobolev inequalities {\juerg to} obtain
\begin{align}
  &{\pier c} \intQt |\phi| \, |\sigmaab - \muab| \, |\sigma|
  \leq {\pier c} \iot \norma\phi_4 \, \norma{\sigmaab - \muab}_4 \, \norma\sigma_2 \, ds
  \non
  \\
  & \leq c \iot \normaV\phi \, \normaV{\sigmaab - \muab} \, \normaH\sigma \, ds
  \leq \delta \iot \normaV\phi^2 \, ds
  + c_\delta \iot \normaV{\sigmaab - \muab}^2 \, \normaH\sigma^2 \, ds
  \non
\end{align}
and we observe at once that the function $s\mapsto\normaV{(\sigmaab-\muab)(s)}^2$
is bounded in $L^1(0,T)$, thanks to \eqref{genest} for $\sigmaab$ 
and \accorpa{genest}{genestbis} for~$\muab$.
Finally, we have {\juerg that}
\Bsist
  && {\pier c} \intQt |\sigma - \mu| \, |\sigma|
  \leq \frac \a 8 \intQt |\mu|^2
  + c  \intQt |\sigma|^2 .
  \non
\Esist
At this point, we collect \eqref{pererroreA} and all the inequalities we have obtained,
rearrange and {\pier infer that
\Bsist
  && \frac 12 \, \normaVp{(\a\mu + \phi + \sigma)(t)}^2
  + \frac \a 8 \intQt |\mu|^2 
  \non
  \\
  && \quad {}
  + \left( 1 -  (1+L)^2 \a -4\delta \right) \iot \normaV{\phi}^2 \, ds  + \frac 12 \iO |\sigma(t)|^2
  + (1-\delta) \iot \normaV\sigma^2 \, ds
  \non
  \\
  \separa
  && \leq  c_\delta \iot \left(1+ \normaV{\sigmaab - \muab}^2 \right) \, \normaH\sigma^2 \, ds
  + c_\delta \iot \normaVp{\a\mu + \phi + \sigma}^2 \, ds 
  +c_\delta \b .
  \non
\Esist
Then, we} choose $\delta$ small enough {\pier
(let us recall that we are assuming $(1+L)^2\a<1$)
and} apply the Gronwall lemma.
This yields \eqref{errorbtozero}, and the proof is complete.{\pier\QED}


\section{Proof of Theorems~\ref{AsymptoticsA} and~\ref{Atozero}}
\label{ATOZERO}
\setcounter{equation}{0}

It is understood that assumptions \HPstruttura\ and \eqref{growth} on the structure 
and \eqref{hpdati} on the initial data are in force.
We first prove Theorem~\ref{AsymptoticsA}.
The main tool is Theorem~\ref{GenEst}, applied to the solution $(\muab,\phiab,\sigmaab,\xiab)$
to problem \Pbl.
However, we need a preliminary estimate {\juerg that} 
has already been performed in~\cite{CGRS}.
Nevertheless, for the reader's convenience, we repeat the core of the argument here.

\step
Auxiliary a priori estimate

We omit the indices $\a$ and~$\b$ for a while.
We observe that the growth condition~\eqref{growth}
formally implies {\pier (see \cite[formulas~(2.25)--(2.26) and Rem.~2.5]{CGRS})}
\Beq
  \iO |\xi(t)|
  \leq c \iO \bigl( \Beta(\phi(t)) + 1 \bigr)
  \quad \aat.
  \label{stimaxi}
\Eeq
Now, we simply integrate \eqref{seconda} over $\Omega$ 
and use the homogeneous Neumann boundary condition for~$\phi$.
We deduce~that
\Beq
  \iO \mu(t) 
  = \iO \b\dt\phi(t) + \iO \xi(t) + \iO \pi(\phi(t))
  \quad \aat .
  \non
\Eeq
As $\pi$ is \Lip\ continuous {\pier and $\Pi$ is nonnegative},
the above identity and \eqref{stimaxi} imply {\juerg that}
\Beq
  \Bigl| \iO \mu(t) \Bigr|
  \leq \iO \b |\dt\phi(t)|
  + c \iO F(\phi(t))
  + c \iO {\pier |\phi(t)|} + c .
  \label{mediamu}
\Eeq
By taking advantage of \eqref{genest},
we deduce that the function $t\mapsto\iO\mu(t)$ is estimated in~$L^2(0,T)$.
Therefore, by {\juerg using} \eqref{genest} once more and the Poincar\'e 
inequality \eqref{poincare},
we conclude that
\Beq
  \norma\mu_{\L2V} \leq c \,.
  \label{damediamu}
\Eeq

\step
Proof of Theorem~\ref{AsymptoticsA}

Now, we {\juerg explicitly} write the indices $\a$ and~$\b$.
First, \eqref{damediamu} ensures that the \lhs s of both \eqref{genest} and \eqref{genestbis} are bounded,
so that, since $\b$ is fixed, we deduce \accorpa{aconvmu}{aconvdt}, at least for a subsequence{\pier :
note that \eqref{aconvdt} follows from \eqref{aconvmu}, \eqref{aconvphi} and the bound for $\{\a \dt \muab + \dt \phiab\}$ in $\L2\Vp$.} 
Now, let $(\mub,\phib,\sigmab,\xib)$ be any limiting quadruplet.
Then, \accorpa{aconvphi}{aconvsigma} imply that the initial conditions for $(\phib,\sigmab)$ are satisfied.
Moreover, the Aubin-Lions lemma (see, e.g., \cite[Thm.~5.1, p.~58]{Lions})
ensures that $\phiab$ converges to $\phib$ strongly in $\L2H$
(even better, of course),
whence, $\pi(\phiab)$ and $p(\phiab)$
converge to $\pi(\phib)$ and $p(\phib)$, respectively, strongly in $\L2H$.
Therefore, we can identify the limits of the nonlinear terms $\xiab$ and~$\Rab$.
For the former, we can apply, e.g., \cite[Cor.~2.4, p.~41]{Barbu}.
For the latter, we note that $\Rab$ converges to $p(\phib)(\sigmab-\mub)$ (at~least) weakly in~$\LQ1$.
At this point, we can write the integrated--in--time version of problem \accorpa{prima}{terza}
for the approximating solution
with time dependent test functions and take the limit as $\a$ tends to zero.
We obtain the analogous systems for $(\mub,\phib,\sigmab,\xib)$,
and this implies \accorpa{primab}{terzab} for such a quadruplet.\QED

Now, we assume that {\pier \accorpa{pconstant}{classico}} hold, in addition,
and start proving Theorem~\ref{Atozero}.
However, as in the previous section, we do not follow the order of the statement.
Indeed, let us assume for a while that its part~$ii)$
has been proved in the following modified version:
the error estimate \eqref{erroratozero} holds for every solution
to the limit problem \Pblb\
(with a constant $C_\b$ that might depend on the solution we are considering).
From this, we derive the uniqueness stated in~$i)$ {\juerg in the following way:
let} $(\mu_i,\phi_i,\sigma_i,\xi_i)$, $i=1,2$, be two solutions.
Then inequality \eqref{erroratozero}
holds for both of them with a common constant~$C_\b$.
Hence, we immediately derive,~e.g.,
\Beq
  \norma{\mu_1-\mu_2}_{\L2V}
  + \norma{\phi_1-\phi_2}_{\L2V}
  + \norma{\sigma_1-\sigma_2}_{\L2V}
  \leq C_\b \, \a^{1/2} 
  \non
\Eeq
for every $\a\in(0,1)$.
We deduce that $\mu_1=\mu_2$, $\phi_1=\phi_2$ and $\sigma_1=\sigma_2$,
whence also $\xi_1=\xi_2$ since $B$ is single-valued.

So, it remains to prove that 
the error estimate \eqref{erroratozero} holds {\juerg true} for every solution
to the limit problem \Pblb\
(with a constant $C_\b$ that might depend on the solution we are considering).
To~this end, we need a further estimate on the solution $(\muab,\phiab,\sigmaab)$ to problem \Pbl.

\step
A new a priori estimate

We proceed formally in order not to be too technical.
Moreover, {\pier let us} omit the indices $\a$ and~$\b$.
We consider equation {\pier \eqref{prima},
as well as the variational equation one obtains by formally differentiating 
\eqref{seconda} with respect to time.}
Then, we present them in the form of abstract equations
by introducing the Riesz operator $\A$ defined by~\eqref{defA},~i.e.,{\pier
\Bsist
  && \dt \bigl( \a \mu + \phi \bigr) 
  + \A \mu 
  = p ( \sigma - \mu) + \mu 
  \non
  \\
  && \dt\mu 
  = \b \dt^2\phi 
  + \A \dt\phi 
  + F''(\phi)\dt\phi 
  - \dt\phi .
  \non
\Esist
}%
{\juerg These} equations are meant in the sense of the Hilbert triplet $(V,H,\Vp)$, \aet.
Now, we test them by $\Am\dt\mu$ {\pier and  $-\Am\dt\phi$, respectively:
we} sum up and integrate over~$(0,t)$.
By omitting the evaluation point inside integrals for brevity, we obtain{\pier
\Bsist
  && \a \iot \< \dt\mu , \Am \dt\mu > \, ds
  + \iot \< \dt\phi , \Am \dt\mu > \, ds
  + \iot \< \A \mu , \Am \dt\mu > \, ds
  \non
  \\
  && \quad {}
  + p \iot \< \mu , \Am \dt\mu > \, ds
  - \iot \< \dt\mu , \Am \dt\phi > \, ds
  + \b \iot \< \dt^2\phi , \Am \dt\phi > \, ds
  \non
  \\
  && \quad {}
  + \iot \< \A \dt\phi , \Am \dt\phi > \, ds
  + \iot \< F''(\phi)\dt\phi , \Am \dt\phi > \, ds
  - \iot \< \dt\phi , \Am \dt\phi > \, ds
  \non
  \\
  \separa
  &&
  = p  \iot \< \sigma , \Am \dt\mu > \, ds
  + \iot \< \mu , \Am \dt\mu > \, ds .
  \non
\Esist
Now, let us recall \accorpa{defA}{propAtre}.
Two terms cancel out by \eqref{propAdue} and we have
\begin{align}
  & \a \iot \normaVp{\dt\mu}^2\, ds
  + \frac 12 \, \normaH{\mu(t)}^2 
   + \frac p2 \, \normaVp{\mu(t)}^2 
  + \frac \b 2 \, \normaVp{\dt\phi(t)}^2
  + \iot \normaH{\dt\phi}^2 \, ds
  \non
  \\
  & = 
  - \iot \bigl( F''(\phi)\dt\phi , \dt\phi \bigr)_* \, ds
  + \iot \normaVp{\dt\phi}^2 \, ds
  + p \iot \bigl( \sigma , \dt\mu \bigr)_* \, ds
  \non
  \\
  & \quad {}
      + \frac 12 \, \normaVp{\mu(t)}^2
      - \frac 12 \, \normaVp\muz^2
    + \frac 12 \, \normaH\muz^2 
    +\frac p2 \, \normaVp\muz^2
  + \frac \b 2 \, \normaVp{\dt\phi(0)}^2\,.
  \label{pier1}
\end{align}
{\juerg Therefore}, we just have to deal with the \rhs.
In order to control the integral involving~$F''$, 
we  account for the dual Sobolev inequality~\eqref{dualsobolev},
the growth condition \eqref{classico}, the \holder\ and Sobolev inequalities,
the estimate \eqref{genest} for~$\phi$, and~have
\Bsist
  && - \iot \bigl( F''(\phi)\dt\phi , \dt\phi \bigr)_* \, ds
  \leq \iot \normaVp{F''(\phi)\dt\phi} \normaVp{\dt\phi} \, ds
  \non
  \\
  && \leq c \iot \norma{F''(\phi)\dt\phi}_{6/5} \, \normaVp{\dt\phi} \, ds
  \leq c \iot \norma{1+\phi^2}_3 \, \norma{\dt\phi}_2 \, \normaVp{\dt\phi} \, ds
  \non
  \\
  && \leq c \iot \bigl( 1 + \norma\phi_6^2 \bigr) \, \norma{\dt\phi}_2 \, \normaVp{\dt\phi} \, ds
  \leq c \iot \bigl( 1 + \normaV\phi^2 \bigr) \, \normaH{\dt\phi} \, \normaVp{\dt\phi} \, ds
  \non
  \\
  && \leq \frac 12 \iot \normaH{\dt\phi}^2 \, ds
  + c \iot \normaVp{\dt\phi}^2 \, ds \,.
  \non
\Esist
Another term to treat is 
\Bsist
  &&  p  \iot \bigl( \sigma , \dt\mu \bigr)_* \, ds = 
  - p \iot \bigl( \dt\sigma , \mu \bigr)_* \, ds 
  + p\, \bigl( \sigma(t) , \mu(t) \bigr)_*
  - p \, \bigl( \sigmaz , \muz \bigr)_*
  \non
  \\
  && \quad{}
  \leq   \frac p4  \iot \normaVp{\mu}^2 \, ds
  + p \, \norma{\dt\sigma}_{\L2\Vp}^2
  + \frac p4  \normaVp{\mu(t)}^2 
  + p \, \norma{\sigma}_{\L\infty\Vp}^2
  + c ,
  \non
\Esist
the last inequality being due to \eqref{damediamu} 
(which holds {\juerg true} under the assumptions of Theorem~\ref{AsymptoticsA},
thus also in the present case) and~\eqref{hpdati}.
The next term on the \rhs\ of \eqref{pier1} is $(1/2)\normaVp{\mu(t)}^2$.
We observe once more that the norm of the embedding $H\subset\Vp$ is~$1$
since we are using the standard norms in $V$ and~$H$.
Therefore, we~have
\Beq
  \frac 12 \, \normaVp{\mu(t)}^2
  \leq \frac 12 \, \normaH{\mu(t)}^2 \,.
  \non
\Eeq
Finally, we consider the norms of the initial values of the time derivatives.
We formally~have
\Beq
  \b \dt\phi(0) = \muz - \A\phiz + \phiz - F'(\phiz)\,,
  \non
\Eeq
which is a fixed element of~$\Vp$ due to~{\pier\eqref{hpdati}}.
Indeed, for~$F'(\phiz)$, we make the following observation.
Assumption \eqref{classico} implies
$|F'(r)|\leq c(|r|^3+1)$ for every $r\in\erre$.
As $\phiz\in V$, we have $\phiz\in\Lx6$ by the Sobolev inequality.
We infer that $\phi^3\in\Lx2$,
whence $F'(\phiz)\in H$.
At this point, we can collect \eqref{pier1} and all the inequalities we have proved, 
then apply the Gronwall lemma. We conclude~that
\Beq \a^{1/2} \norma{\dt\muab}_{\L2\Vp} + p^{1/2} \norma{\muab}_{\L\infty\Vp}
  + \norma{\dt\phiab}_{\L\infty\Vp}
  \leq c \,,
  \label{newestimate}
\Eeq
where we have used the full notation with indices.}

\step
Conclusion of the proof of Theorem~\ref{Atozero}

We prove that \eqref{erroratozero} holds for every solution $(\mub,\phib,\sigmab)$ to problem \Pblb\
(with a constant $C_\b$ that might depend on the solution we are considering).
We often omit writing the evaluation point, {\pier explicitly}, in order to simplify the notation.
We take the difference between \accorpa{prima}{terza}, written for 
$(\muab,\phiab,\sigmaab)$,
and \accorpa{primab}{terzab}, written for $(\mub,\phib,\sigmab)$, 
where {\pier $\xiab=B(\phiab)$ and {\gianni $\xib=B(\phib)$}.}
By setting, for convenience,
\Beq
  \mu := \muab - \mub \,, \quad
  \phi := \phiab - \phib
  \aand
  \sigma := \sigmaab - \sigmab\,, 
  \non
\Eeq
and recalling that $p$ is a {\pier nonnegative constant (cf.~\eqref{pconstant})},
we~have
\Bsist
  \hskip-1cm
  && \iO \dt\phi \, v + \iO \nabla\mu \cdot \nabla v
  = - \a \< \dt\muab , v >
  + p \iO (\sigma - \mu) v 
  \non
  \\
  && \iO \mu v
  = \b \iO \dt\phi \, v
  + \iO \nabla\phi \cdot \nabla v
  + \iO (\xiab - \xib) v
  + \iO \bigl(\pi(\phiab) - \pi(\phib)) v
  \non
  \\
  && \< \dt\sigma , v > + \iO \nabla\sigma \cdot \nabla v
  = p \iO (\mu - \sigma) v \,,
  \non
\Esist
where each {\pier equality} holds {\juerg true} for every $v\in V$ and \aet.
Now, we {\pier take the test function $v$ equal to}
\Beq
  \phi + \b\mu, \quad
  \mu - \phi
  \aand
  \sigma\,,
  \non
\Eeq
respectively.
Next, we integrate with respect to time, {\pier exploit} the initial conditions
$\phi(0)=0$ and $\sigma(0)=0$,
add the volume integrals of $p\b|\mu|^2$ and of $p|\sigma|^2$ to both sides for convenience, 
and rearrange a little. 
We have, for every $t\in[0,T]$,
\Bsist
  && \frac 12 {\pier \iO } |\phi(t)|^2
  + \b \intQt \dt\phi \, \mu
  + \intQt \nabla\mu \cdot \nabla\phi
  + \b \intQt |\nabla\mu|^2
  \non
  \\
  && \quad {}
  + (1 + p\b) \intQt |\mu|^2
  + \frac \b 2{\pier \iO} |\phi(t)|^2
  - \b \intQt \dt\phi \, \mu
  + \intQt |\nabla\phi|^2
  - \intQt \nabla\phi \cdot \nabla\mu
  \non
  \\
  && \quad {}
  + \intQt (\xiab - \xib) \phi
  + \frac 12 \iO |\sigma(t)|^2
  + \intQt |\nabla\sigma|^2
  + p \intQt |\sigma|^2
  \non
  \\
  \separa
  && = \intQt \mu \phi
  - \intQt \bigl( \pi(\phiab) - \pi(\phib) \bigr) \phi
  + \intQt \bigl( F'(\phiab) - F'(\phib) \bigr) \mu
  \non
  \\
  && \quad {}
  - \a \iot \< \dt\muab , \phi + \b\mu > \, ds
  + p \intQt \bigl\{ (\sigma - \mu) (\phi + \b\mu - \sigma) + \b|\mu|^2 + |\sigma|^2 \bigr\} .
  \qquad
  \label{pererrorea}
\Esist
Four terms on the \lhs\ cancel out, and the other ones are nonnegative.
Now, we estimate each integral on the \rhs, separately.
For the first two of them, we~have
\Beq
  \intQt \mu \phi
  - \intQt \bigl( \pi(\phiab) - \pi(\phib) \bigr) \phi
  \leq \frac 14 \intQt |\mu|^2 
  + c \intQt |\phi|^2 
  \non
\Eeq
since $\pi$ is \Lip\ continuous.
In order to deal with the next integral, we observe that
\Beq
  |F'(\phiab) - F'(\phib)|
  \leq \bigl( |\phiab|^2 + |\phib|^2 + 1 \bigr) |\phi|
  \quad \aeQ\,,
  \non
\Eeq
due to the mean value theorem and~\eqref{classico}.
Hence, applying the \holder\ and Sobolev inequalities,
and owing to estimate~\eqref{genest} for $\phiab$ and the regularity \eqref{regphib} of~$\phib$,
we {\juerg can infer that} 
\Bsist
  && \intQt \bigl( F'(\phiab) - F'(\phib) \bigr) \mu
  \leq c \iot \bigl\| |\phiab|^2 + |\phib|^2 + 1 \bigr\|_3 \, \norma\phi_2 \, \norma\mu_6 \, ds
  \non
  \\
  && \leq c \iot \bigl( \norma\phiab_6^2 + \norma\phib_6^2 + 1 \bigr) \, \norma\phi_2 \, \norma\mu_6 \, ds
  \non
  \\
  && \leq c \iot \normaH\phi \, \normaV\mu \, ds  
  \leq \delta \iot \normaV\mu^2 \, ds
  + c_\delta \intQt |\phi|^2 \,.
  \non
\Esist
Next, we take advantage of \eqref{newestimate} for $\dt\muab$
and deduce that
\Bsist
  && - \a \iot \< \dt\muab , \phi + \b\mu > \, ds
  \leq \delta \iot \bigl( \normaV\phi^2 + \normaV\mu^2 \bigr) \, ds
  + c_\delta \, \a^2 \iot \normaVp{\dt\muab}^2 \, ds
  \non
  \\
  && \leq \delta \iot \bigl( \normaV\phi^2 + \normaV\mu^2 \bigr) \, ds
  + c_\delta \, \a \,.
  \non
\Esist
Finally, the last {\gianni term} {\pier in \eqref{pererrorea}} can be simplified and estimated 
{\juerg in the following way:}
\Bsist
   &&{\pier  p \intQt \bigl\{ (\sigma - \mu) (\phi + \b\mu - \sigma) + \b|\mu|^2 + |\sigma|^2 \bigr\}}
   \non
   \\
  &&
  \leq {\pier p }\intQt \bigl( \sigma\phi - \mu\phi + (\b+1) \sigma\mu \bigr)
  \leq \delta \intQt |\mu|^2
  + c_\delta \intQt \bigl( |\phi|^2 + |\sigma|^2 \bigr) .
  \non
\Esist
By combining \eqref{pererrorea} and these inequalities, 
choosing $\delta$ small enough and applying the Gronwall lemma,
we obtain~\eqref{erroratozero}.



\vspace{3truemm}

\Begin{thebibliography}{10}

\bibitem{Barbu}
{\sc V. Barbu},
``Nonlinear Differential Equations of Monotone Types in Banach spaces'',
Springer Monographs in Mathematics, 2010.

\bibitem{Brezis}
{\sc H. Brezis,}
``Op\'erateurs maximaux monotones et semi-groupes de contractions
dans les espaces de Hilbert'',
North-Holland Math. Stud.~{\bf 5},
North-Holland, Amsterdam, 1973.

\bibitem{CH} 
{\sc J.W.~Cahn and J.E.~Hilliard}, 
{\em Free energy of a nonuniform system I. Interfacial free energy}, 
J. Chem. Phys.,
{\bf 2} (1958), pp.~258--267.

{\pier
\bibitem{CC1}
{\sc G. Canevari and P. Colli}, 
{\em Solvability and asymptotic analysis
of a generalization of the Caginalp phase field system}, 
Commun. Pure Appl. Anal., {\bf 11} (2012), pp.~1959--1982.

\bibitem{CC2}
{\sc G. Canevari and P. Colli}, 
{\em Convergence properties for a generalization of the Caginalp phase field system}, 
Asymptot. Anal., {\bf 82} (2013), pp.~139--162.
} 
 
{\pier
\bibitem{CGG}
{\sc P. Colli, G. Gilardi, and M. Grasselli},
{\em Asymptotic analysis of a phase field model with memory for vanishing time
relaxation}, 
Hiroshima Math. J., {\bf 29} (1999), pp.~117--143.
}

{\pier
\bibitem{CGH}
{\sc P. Colli, G. Gilardi, and D. Hilhorst},
{\em On a Cahn--Hilliard type phase field system related to tumor growth},
{\pier Discrete Contin. Dyn. Syst.,} \textbf{35} (2015), pp.~2423--2442.}

{\pier
\bibitem{CGPS}
{\sc P. Colli, G. Gilardi, P. Podio-Guidugli, and J. Sprekels,}
{\em An asymptotic analysis for a nonstandard Cahn-Hilliard system with viscosity},
Discrete Contin. Dyn. Syst. Ser. S, \textbf{6} (2013), pp.~353--368.
}

\bibitem{CGRS}
{\sc P. Colli, G. Gilardi, E. Rocca, {\pier and} J. Sprekels},
{\pier {\em Vanishing viscosities and error estimate 
for a Cahn--Hilliard type phase field system 
related to tumor growth},
preprint arXiv:1501.07057 [math.AP] (2015) 1--19.} 
 
{\pier
\bibitem{CS}
{\sc P. Colli and J. Sprekels},
{\em Stefan problems and the
Penrose--Fife phase field model}, 
Adv.  Math. Sci. Appl., {\bf 7} (1997), pp.~911--934.
}
 
{\pier
\bibitem{DKS}
{\sc A. Damlamian, N. Kenmochi, and N. Sato},
{\em Subdifferential operator approach to a class 
of nonlinear systems for Stefan problems with phase relaxation},
Nonlinear Anal., {\bf 23} (1994), pp.~115--142.
}

\bibitem{EllSt} 
{\sc C.M. Elliott and A.M. Stuart}, 
{\em Viscous Cahn--Hilliard equation. II. Analysis}, 
J. Differential Equations, {\bf 128} (1996),  pp.~387--414.

\bibitem{EllZh} 
{\sc C.M. Elliott and S. Zheng}, 
{\em On the Cahn--Hilliard equation}, 
Arch. Rational Mech. Anal., {\bf 96} (1986), pp.~339--357.

{\pier
\bibitem{FGR}
{\sc S. Frigeri, M. Grasselli, and E. Rocca},
{\em On a diffuse interface model of tumor growth},
European J. Appl. Math., DOI: 10.1017/S0956792514000436
}
 
{\pier
\bibitem{Gir}
{\sc M. Girotti},
{\em Vanishing time relaxation for a phase-field 
model with entropy balance},
Adv. Math. Sci. Appl., {\bf 22} (2012), pp.~553--575.
}

{\betti 
\bibitem{HZO} 
{\sc A. Hawkins-Daarud, K.G. van der Zee, and J.T. Oden}, 
{\em Numerical simulation of a thermodynamically consistent 
four-species tumor growth model}, Int. J. Numer. Meth. Biomed. 
Engng., \textbf{28} (2011), pp.~3--24.
}

{\pier
\bibitem{HKNV} 
{\sc D. Hilhorst, J. Kampmann, T.N. Nguyen, and K.G. Van der Zee}, 
{\em Formal asymptotic limit of a diffuse-interface tumor-growth model},
Math. Models Methods Appl. Sci.,  DOI: 10.1142/S0218202515500268
}

\bibitem{Lions}
{\sc J.-L.~Lions},
``Quelques m\'ethodes de r\'esolution des probl\`emes
aux limites non lin\'eaires'',
Dunod; Gauthier-Villars, Paris, 1969.
 
{\pier
\bibitem{Roc}
{\sc E. Rocca},
{\em Asymptotic analysis of a conserved phase-field model with memory for vanishing time relaxation}, 
Adv. Math. Sci. Appl., {\bf 10} (2000), pp.~899--916. 
}

{\pier
\bibitem{Ros1}
{\sc R. Rossi},
{\em Asymptotic analysis of the Caginalp phase-field model 
for two vanishing time relaxation parameters}, 
Adv. Math. Sci. Appl., \textbf{13} (2003), pp.~249--271. 
}

{\pier
\bibitem{Ros2}
{\sc R. Rossi},
{\em Well-posedness and asymptotic analysis for a Penrose-Fife type phase field system},
Math. Methods Appl. Sci., \textbf{27} (2004), pp.~1411--1445. 
}
 
{\pier
\bibitem{Sch}
{\sc G. Schimperna},
{\em Singular limit of a transmission problem for the parabolic phase-field model}, Appl. Math., {\bf 45} (2000), pp.~217--238.
}

{\betti 
\bibitem{WZZ} 
{\sc X. Wu, G.J. van Zwieten, and K.G. van der Zee}, 
{\em Stabilized second-order convex splitting schemes for Cahn--Hilliard
models with applications to diffuse-interface tumor-growth models}, 
Int. J. Numer. Meth. Biomed. Engng., \textbf{30} (2014), pp.~180--203.
}

\End{thebibliography}

\End{document}

\bye